\documentclass{article}
 \newtheorem{thm}{Theorem}[section]

\usepackage{amssymb}
\begin{document}
\title{Analytic Classification of Plane Branches up to Multiplicity 4}

\author{\begin{tabular}{ll}
{\normalsize Abramo Hefez} \thanks{partially supported by CNPq,
PRONEX and PROCAD}   &
{\normalsize Marcelo E. Hernandes} \thanks{partially supported by CAPES and PROCAD} \\
{\normalsize Universidade Federal Fluminense} & {\normalsize Universidade Estadual de Maring\'a} \\
{\normalsize Instituto de Matem\'atica} & {\normalsize Departamento de Matem\'atica} \\
{\normalsize R. Mario Santos Braga, s/n } & {\normalsize Av. Colombo, 5790} \\
{\normalsize 24020-140 Niter\'oi, RJ - Brazil} & {\normalsize 87020-020 Maring\'a, PR - Brazil} \\
{\small {\em E-mail}:\quad hefez@mat.uff.br} & {\small {\em
E-mail}:\quad mehernandes@uem.br}
\end{tabular}}
\pagestyle{myheadings} \markboth{HEFEZ AND
HERNANDES}{Classification of Plane Branches up to Multiplicity 4}

\date{ \ }

\maketitle

\begin{abstract}
We perform the analytic classification of plane branches of
multiplicity less or equal than four. This is achieved by
computing a Standard basis for the modules of K\"ahler
differentials of such branches by means of the algorithm we
developed in \cite{[HH1]} and then applying the classification
method for plane branches presented in \cite{[HH2]}.\vspace{3mm}



\end{abstract}


\section{Introduction}

Until very recently, the analytic classification of plane branches
within an equisingularity class was an open problem (cf.
\cite{[HH2]}). The first serious attempt to solve this problem was
made by S. Ebey in \cite{[E]}, where he classified branches of
multiplicity two and three and very few classes of branches of
multiplicity four. Few years later, O. Zariski
dedicated the book \cite{[Z2]} to the study of this problem
without much success.

In this paper we will exploit the general approach of
\cite{[HH2]}, where a general method to perform the effective
analytic classification of plane branches within a given
equisingularity class is described, to classify all branches with
multiplicity less or equal than four. An important invariant in
this context, used to stratify the parameter space $\Sigma$ of any
given equisingularity class, is the value set $\Lambda$ of the
associated module of K\"ahler differentials, which we determine
computing a Standard Basis for the module by means of the
algorithm we developed in \cite{[HH1]}.

Some authors have used other invariants as Tjurina's number $\tau$
of the branch (cf. \cite{[LP]}) or Zariski's $\lambda$ invariant
(cf. \cite{[P]}), introduced in \cite{[Z1]}, to stratify the
parameter space $\Sigma$. These invariants express very partial
information one can get from the set $\Lambda$ and do not separate
properly branches for classification purpose. Another invariant
adopted in \cite{[GP]} and \cite{[GHP]} is the Hilbert function on
the Tjurina algebra of the branch. This also doesn't work properly
and in Section 6 we will give an example in which two branches
with two different sets $\Lambda$, have the same Hilbert function.

It should be noted that the differentials obtained by means of the
above mentioned algorithm may also serve to determine explicitly
the analytic coordinate changes that will reduce a given curve to
its normal form (cf. \cite{[HH2]} and \cite{[Her]}). It also
should be noted that our classification corresponds to the
classification by contact equivalence, introduced by J. Mather,
and is not the same as Arnold's (see \cite{[Ar]} and \cite{[J]}),
where germs of functions are classified up to right equivalence;
i.e., up to changes of coordinates in the source. It is worth
noting that the contact equivalence classification is much more
difficult than the right equivalence classification.

Finally, we should mention that in \cite{[NW]} the Tjurina's
numbers were computed for all plane branches of multiplicity less
or equal than four. The method used there was to determine, by
ad-hoc calculations, the cardinality of the set $\Lambda\setminus
\Gamma$, where $\Gamma$ is the semigroup of values of the branch,
which measures the difference between the conductor of $\Gamma$
and Tjurina's number (cf. Remark 4.8 of \cite{[HH1]}). This can
easily be derived from our results since we give explicitly, for
such branches, the sets $\Lambda$ and $\Gamma$.

\section{Preliminaries}

Let $A$ be either the ring of formal or convergent power series in
two indeterminates $X$ and $Y$ with coefficients in $\mathbb C$
and $B$ the ring of formal or convergent power series in one
indeterminate $t$ with coefficients in $\mathbb C$. A plane branch
$(f)$ is a class in $A$, modulo associates, of an irreducible
non-unit $f$ in $A$. Two branches $(f_1)$ and $(f_2)$ are {\em
analytically equivalent}, or shortly {\em equivalent}, if there
exist an automorphism $\Phi$ and a unit $u$ of $A$ such that
$\Phi(f_1)=uf_2$.

Consider a Puiseux parametrization $\varphi=(x(t),y(t))\in B\times
B$ of the branch $(f)$, and its associated map germ
$\varphi:(\mathbb C,0) \rightarrow (\mathbb C^2,0)$. It is known
(see Lemma 2.2 in \cite{[BG]}) that given two plane branches
$(f_1)$ and $(f_2)$ parame\-tri\-zed, respectively, by $\varphi_1$
and $\varphi_2$, then $(f_1)$ and $(f_2)$ are analytically
equivalent if, and only if, $\varphi_1$ and $\varphi_2$ are
${\mathcal A}$-equivalent, where $\mathcal A$-equivalence means
that there exist germs of analytic isomorphisms $\sigma$ and
$\rho$ of $(\mathbb C^2,0)$ and $(\mathbb C,0)$, respectively,
such that $\varphi_2=\sigma\circ \varphi_1 \circ \rho^{-1}$.

So, the analytic classification of plane branches reduces to the
$\mathcal A$-classification of parametrizations, which we are
going to undertake in this paper.

To the map germ $\varphi:(\mathbb C,0) \rightarrow (\mathbb
C^2,0)$ there is associated a ring homomorphism
$\varphi^*:A\rightarrow B$, determining a natural valuation
$v_\varphi$ on $A$. The value set $\Gamma=v_\varphi(A)\subset
\mathbb N \cup \{\infty\}$ will be called the {\em semigroup of
values} of the branch. This is a well known complete invariant for
the topological classification of plane branches.

The semigroup $\Gamma$ of a plane branch has a conductor $c$ and
any element in the finite set $\mathbb N \setminus \Gamma$ is
called a gap of $\Gamma$. If $v_0< v_1< \cdots < v_g$ is a minimal
set of generators for $\Gamma$, then the multiplicity of the branch is
$v_0$ and one has the formula (cf. \cite{[H]}, (7.1)):
\begin{equation} c=\sum_{i=1}^g(n_i-1)v_i-v_0+1,
\end{equation}
where $n_i=\displaystyle \frac{e_{i-1}}{e_i}$, and $e_i={\rm
GCD}(v_0,\ldots,v_{i})$, $1\leq i \leq g$.

There is also an induced modules homomorphism $\varphi^*:AdX+AdY
\rightarrow Bdt$,
$\varphi^*(gdX+hdY)=(g(x(t),y(t))x'(t)+h(x(t),y(t))y'(t))dt$,
which induces a valuation, also denoted by $v_\varphi$, defined by
$$v_\varphi(gdX+hdY)=v_\varphi(g(x(t),y(t))x'(t)+h(x(t),y(t))y'(t))+1.
$$

This allows to define the value set $\Lambda=v_\varphi(AdX+AdY)$,
which is an analytic invariant of the branch (see \cite{[HH2]},
Proposition 3.2). Notice that $\Gamma\subset \Lambda\cup \{0\}$.

Zariski in \cite{[Z1]} has shown that $\Lambda\setminus
\Gamma=\emptyset$ if, and only if, the branch is analytically
equivalent to one with Puiseux parametrization
$(t^{v_0},t^{v_1})$, with GCD$(v_0,v_1)=1$.

On the other hand, if $\Lambda\setminus \Gamma \neq \emptyset$, Zariski also
showed that the branch is analytically equivalent to one with
Puiseux parametrization
\[
(t^{v_0},t^{v_1}+t^\lambda+\sum_{\lambda<i<c-1}a_it^i),
\]
where $\lambda=\min(\Lambda\setminus \Gamma)-v_0$, is the so
called {\em Zariski's invariant}. The invariant $\lambda$ may take
any value in the following set
\[
\{ i\in \mathbb N; v_1<i< v_2-(n_1-1)v_1, \ e_1|i \ \} \cup \{
v_2-(n_1-1)v_1\}.
\]

In \cite{[HH2]} we proved the following stronger
result:\vspace{2mm}

\noindent{\sc The Normal Forms Theorem} {\bf (NFT)} \quad {\em Let $(f)$ be a
plane branch with semigroup of values $\Gamma=\langle
v_0,v_1,\ldots,v_g\rangle$ and value set of differentials
$\Lambda$. If $\Lambda\setminus\Gamma\neq \emptyset$, then $(f)$
is equivalent to a branch with a Puiseux parametrization
\begin{equation}
(t^{v_0},t^{v_1}+t^\lambda+\sum_{\stackrel{i>\lambda}{i\not \in
\Lambda-v_0}} a_it^i).
\end{equation}
Moreover, if $\varphi$ and $\varphi'$ {\rm (}with coefficients
$a_i'$ instead of $a_i${\rm )} are parametrizations as in {\rm
(2.2)}, representing two plane branches $(f)$ and $(g)$ with same
semigroup of values and same set of values of differentials, then
$(f)$ is equivalent to $(g)$ if, and only if, there exists $r\in
\mathbb C^*$ such that $r^{\lambda-v_1}=1$ and
$a_i=r^{i-v_1}a_i'$, for all $i$.} \vspace{2mm}

In \cite{[HH1]} we developed an algorithm to compute a Standard
basis for the $\varphi^*(A)$-module $\varphi^*(AdX+AdY)$, whose
values added to the elements of $\Gamma$ give the set $\Lambda$.
So, this allows to determine all possible sets $\Lambda$ for the
plane branches with given semigroup $\Gamma$, making effective
the above theorem.

\section{Branches of Multiplicity Less than Four}

The case of multiplicity one will be disregarded since all such
branches are equivalent to each other (cf. \cite{[H]} Proposition 3.1).

Now suppose that a plane branch with multiplicity $v_0=2$ is
given. Then, its semigroup of values is given by $\Gamma=\langle
2, v_1\rangle$, with $v_1$ odd. According to formula (2.1), the
conductor of $\Gamma$ is $c=v_1-1$. So,
$\Lambda\setminus\Gamma=\emptyset$ and by Zariski's result we have
that the given branch is equivalent to one with Puiseux
parametrization $(t^2, t^{v_1})$.

This gives the classification of all branches of multiplicity two.
\vspace{2mm}

Let a branch of multiplicity $v_0=3$ be given, then in this case,
$\Gamma=\langle 3, v_1 \rangle$, with GCD$(3,v_1)=1$, whose
conductor is $c=2(v_1-1)$. The gaps of $\Gamma$ above $v_1$ are
the numbers:
\[
 2v_1-3\left[ \frac{v_1}{3}\right], 2v_1-3\left( \left[
\frac{v_1}{3}\right] -1\right), \ldots, 2v_1-3\cdot 2, 2v_1-3.\]

If $\Lambda\setminus \Gamma=\emptyset$, then the branch is
equivalent to one with a parametrization $(t^{3},t^{v_1})$.

If $\Lambda\setminus \Gamma\neq \emptyset$, then the invariant
$\lambda$ may be any of the following integers:
\[ 2v_1-3\left[ \frac{v_1}{3}\right], 2v_1-3\left( \left[
\frac{v_1}{3}\right] -1\right), \ldots, 2v_1-3\cdot 2 .\]

Once $\lambda$ is chosen, in the above set, it follows that any
gap $j>\lambda$, is such that $j\in \Lambda$. Hence, by the {\bf
NFT} (Normal Forms Theorem) we have that the given branch is
equivalent to one with Puiseux parametrization
\[
(t^{3},t^{v_1}+t^\lambda).
\]

Clearly, two parametrizations as above are equivalent if, and only
if, they are identical.

To describe the set $\Lambda\setminus \Gamma$, suppose, for
example, that $\lambda=2v_1-3k$, for $2\leq k \leq
\left[\frac{v_1}{3}\right]$, then
$$\Lambda\setminus\Gamma=
\left\{ 2v_1-3j; \ 1\leq j \leq k-1 \right\}.$$

\section{Branches of Multiplicity 4}

We will describe in this section all possible sets $\Lambda$ for a
plane branch of multiplicity $v_0=4$. It is easy to verify that
the only possible semigroups of values of multiplicity 4 are
either of the form $\langle 4,v_1\rangle$, with GCD$(4,v_1)=1$,
or of the form $\langle 4,v_1,v_2\rangle$, with
GCD$(4,v_1)=2$ and GCD$(4,v_1,v_2)=1$.

We will consider first the case of semigroups of the form $\Gamma =\langle
4,v_1\rangle$. According to (2.1), we have $c=3(v_1-1)$.

Assume that $\Lambda\setminus \Gamma\neq \emptyset$;
otherwise, the branch would be equivalent to one with
Puiseux parametrization $(t^4, t^{v_1})$.

The gaps of $\Gamma$, above $v_1$, are of the form
$$
2v_1-4j,  \  1\leq j \leq \left[ \frac{v_1}{4} \right];  \ \
3v_1-4j,   \ 1\leq j \leq \left[ \frac{v_1}{2} \right].
$$

So, the $\lambda$ invariant may take any of the following values:
$$
2v_1-4j,  \  2\leq j \leq \left[ \frac{v_1}{4} \right];  \ \
3v_1-4j,   \ 2\leq j \leq \left[ \frac{v_1}{2} \right].
$$

We will analyze two distinct cases according to the value
of $\lambda$:

\noindent {\bf Case a)} $\lambda =3v_1-4j$, for some $j=2,\ldots
,\left [ \frac{v_1}{2}\right ]$.

Using that $\lambda+v_0=\min(\Lambda\setminus \Gamma)$ and the
{\bf NFT}, we may assume that the branch is equivalent to one with
the following Puiseux parametrization:
$$
\varphi(t)=\left(t^4, t^{v_1}+t^{3v_1-4j}+\sum_{i=1}^{j-\left [
\frac{v_1}{4}\right ]-2}{a_it^{2v_1-4\left ( j-\left [
\frac{v_1}{4}\right ] -i\right )}}\right),
$$
because $\{3v_1-4k; \ 1\leq k \leq j-1\}\subset \Lambda$.

 Applying now to these branches Algorithm 4.10 of
\cite{[HH1]}, which will be referred in the sequel as {\em the
algorithm}, simply, we get only one minimal non-exact differential
(MNED) $\omega_1 =xdy-\frac{v_1}{4}ydx$ with
$v_\varphi(\omega_1)=3v_1-4(j-1)$.

Hence,
$$\Lambda\setminus\Gamma =\{ v_\varphi(\omega_1)+\gamma \not \in
\Gamma; \ \gamma \in \Gamma \} = \{3v_1-4s; \ 1 \leq s \leq
j-1\}.$$

\noindent {\bf Case b)} $\lambda =2v_1-4j$, for some $j=2,\ldots
,\left [ \frac{v_1}{4}\right ]$.

Applying the {\bf NFT}, remembering that $\min(\Lambda\setminus
\Gamma)=\lambda+v_0$, we may assume that the branch is equivalent
to one with the following Puiseux parametrization:
\[
\varphi(t)=\left(t^4,t^{v_1}+t^{2v_1-4j}+\sum_{i=1}^{\left [
\frac{v_1}{4}\right ]}{a_it^{3v_1-4\left ( \left [
\frac{v_1}{4}\right ] +j+1-i\right )}}\right),
\]
because $\{ 2v_1-4k, 3v_1-4k; \ 1\leq k \leq j-1\}\subset
\Lambda$.

Applying the algorithm to the above parametrization, we get in the
first step the MNED $\omega_1=xdy-\frac{v_1}{4}ydx$, for
which\vspace{1mm}

$\displaystyle\frac{\varphi^* (\omega_1)}{dt}=(v_1-4j)t^{2v_1-4(j-1)-1}$

$\hspace{1.5cm} + \sum_{i=1}^{\left [ \frac{v_1}{4}\right ]}{\left
( 2v_1-4\left ( \left [ \frac{v_1}{4}\right ] +j+1-i\right )\right
)a_it^{3v_1-4\left ( \left [ \frac{v_1}{4}\right ] +j-i\right
)-1}}.$\vspace{1mm}

In the second step of the algorithm, we have
$$\omega_2=v_1x^{j-1}\omega_1-(v_1-4j)ydy.$$

Hence,
$$\displaystyle\frac{\varphi^*(\omega_2)}{dt}=\sum_{i=1}^{\left [
\frac{v_1}{4}\right ]-j}{v_1a_i\left ( 2v_1-4\left ( \left [
\frac{v_1}{4}\right ] +j+1-i\right )\right )t^{3v_1-4\left ( \left
[ \frac{v_1}{4}\right ] +1-i\right ) -1}} $$
$$
\hspace*{2.5cm}+ \left ( v_1a_{\left [ \frac{v_1}{4}\right ]+1-j}
\left ( 2v_1-8j\right ) -\left ( v_1-4j\right )\left (
3v_1-4j\right )\right )t^{3v_1-4j-1}+\cdots $$

If, for some $i=1,\ldots ,\left [ \frac{v_1}{4}\right ]-j$, we
have $a_i\not =0$, then $\omega_2$ is a MNED, with
$v_\varphi(\omega_2)=3v_1-4\left ( \left [ \frac{v_1}{4}\right ]
+1-k\right )$, where $k=\mbox{min}\{ i;\ a_i\not =0\}$. Hence, the
algorithm ends since we already got $v_0-2=2$ MNED's (cf. Remark
4.11 of \cite{[HH1]}).

In this case, we have that

$\Lambda\setminus\Gamma =\{
v_\varphi(\omega_1)+\gamma\not\in\Gamma ;\ \gamma\in\Gamma
\}\cup\{ v_\varphi(\omega_2)+\gamma\not\in\Gamma ;\
\gamma\in\Gamma\}=$

$\hspace{1.25cm} \{ 2v_1-4s;\ 1\leq s\leq j-1\}\cup\{ 3v_1-4s;\
1\leq s\leq\left [ \frac{v_1}{4}\right ] +1-k \}.$

If $a_i=0$ for all $i$ with $1\leq i\leq \left [
\frac{v_1}{4}\right ]-j$ and $a_{\left [ \frac{v_1}{4}\right
]+1-j}\not =\frac{v_1+\lambda }{2v_1}$, then $\omega_2$ is a MNED
with $v_\varphi(\omega_2)=3v_1-4j$. Again, the algorithm ends since we
already got the maximum number of MNED's.

In this situation,
$$\Lambda\setminus\Gamma =\{ 2v_1-4s;\ 1\leq s\leq j-1\}\cup\{ 3v_1-4s;\ 1\leq s\leq j \}.$$

On the other hand, if $a_i=0$ for all $i=1,\ldots ,\left [
\frac{v_1}{4}\right ] -j$ and $a_{\left [ \frac{v_1}{4}\right
]+1-j}=\frac{v_1+\lambda }{2v_1}$, then the algorithm ends. In
this case, $\omega_1$ is the only MNED, and
$$\Lambda\setminus\Gamma =\{ 2v_1-4s;\ 1\leq s\leq j-1\}\cup\{ 3v_1-4s;\ 1\leq s\leq j-1 \}.$$

We will analyze now the case $\Gamma =\langle 4,v_1,v_2\rangle$.
It is easy to see that $\lambda =v_2+v_1-4m_1$, where
$m_1=\frac{v_1}{2}$ that is, $\lambda=v_2-v_1$, and that the gaps
of $\Gamma$ above $v_1$ are the following:
$$v_2-4j,  1\leq j \leq \left[ \frac{v_2-v_1}{4} \right];  \ \
v_2+v_1-4j, 1\leq j \leq \left[ \frac{v_2}{4} \right].
$$

By the {\bf NFT}, a branch with the above $\Gamma$ is equivalent
to one with Puiseux parametrization
\[
\varphi(t)=(t^4,t^{v_1}+t^{\lambda}+\sum_{i=1}^{\left [
\frac{v_1}{4}\right ]-1}{a_it^{v_2-4\left ( \left [
\frac{v_1}{4}\right ]+1-i\right )}}), \] because $\{v_2+v_1-4k;
1\leq k \leq m_1-1\}\subset \Lambda$.

Putting $z=y^2-x^{m_1}$, we have that $v_\varphi(z)=v_2$ and
therefore $\{x,y,z\}$ is a minimal Standard basis for $\varphi^*(A)$.

In the first step of the algorithm we get
$\omega_1=xdy-\frac{v_1}{4}ydx$ and
$\omega_2=2ydy-m_1x^{m_1-1}dx$. But, $\omega_2=dz$, hence it is an
exact differential. On the other hand, $\omega_1$ is a MNED,
namely $$
\begin{array}{ll}
\displaystyle\frac{\varphi^* (\omega_1)}{dt} = &
 (v_2-4m_1)t^{v_2+v_1-4(m_1-1)-1}
\\ & \\
&
 + \sum_{i=1}^{\left [ \frac{v_1}{4}\right ]-1} \left (
v_2-v_1-4\left ( \left [ \frac{v_1}{4}\right ]+1-i\right )\right
)a_it^{v_2-4\left ( \left [ \frac{v_1}{4}\right ]-i\right )-1}
\end{array}
$$ with $v_\varphi(\omega_1)=v_2+v_1-4(m_1-1)$. Moreover, $\omega_1$ is the
unique MNED, and we have
$$\Lambda\setminus\Gamma =\{ v_2+v_1-4s;\ 1\leq s\leq m_1-1\}.$$

Now, a further application of the {\bf NFT} allows to deduce the following theorem

\begin{thm} A plane branch of multiplicity less or equal than $4$ is
equivalent to a member $C_a$ of one of the families described in
Table {\rm (1)}. Two branches $C_a$ and $C_{a'}$ belonging to
distinct families are never equivalent, and if they belong to the
same family, they are equivalent if and only if they differ by an
homothety; that is, there is $r\in \mathbb C^*$ with
$r^{\lambda-v_1}=1$, where $\lambda$ is Zariski's invariant of the
branch, such that $a'_i=r^{i-v_1}a_i$, for all $i>\lambda$.
\end{thm}

\section{Examples and Remarks}

The method of classification of branches we used in this paper is
effective; that is, it is possible to perform the computations in
order to put any given branch into normal form. For a computer
implementation in MAPLE see \cite{[Her]}.

In \cite{[Z2]}, Section 3, Chapter IV, after a long computation,
Zariski concluded that all branches in the equisingularity class
determined by the semigroup $\langle 4, 6,v_2\rangle$ are
equivalent to each other. This follows immediately from the last
row in Table (1), from where we see that all such branches are
equivalent to the branch with Puiseux parametrization
$$ \varphi(t)=(t^4,t^6+t^{v_2-6}).
$$

In \cite{[E]}, Ebey, and subsequently P. Carbonne in \cite{[Ca]}
studied the equisingularity classes determined by $\langle 4, v_1
\rangle$, where $v_1\leq 11$. All their results are contained in
Table (1).

In \cite{[BH]}, the family of
branches
$$ \varphi_a(t)=(t^4,t^{v_1}+t^{2v_1-8}+\frac{3v_1-8}{2v_1}t^{3v_1-16}+at^{3v_1-12}),$$
was considered. At that time, we hadn't the tools to classify the
members of this family, modulo analytic equivalence. Now, since we
are in the last case of the semigroup $\langle 4, v_1 \rangle$ in
Table (1), we have that $\varphi_a$ determines a plane branch
equivalent to that determined by
$\varphi_{a'}$ if and
only if $a'= r^{2v_1-12}a$, where $r^{v_1-8}=1$.

As mentioned in the introduction, we will now show that $\Lambda$,
the set of values of K\"ahler differentials of the local ring of a
plane branch, is a finer
invariant than the Hilbert function on the Tjurina algebra of that
local ring, used in \cite{[GP]} and \cite{[GHP]}.

To illustrate this, consider the branches given by
\[
f(X,Y)=Y^4-X^9+X^7Y, \ \ \ \ g(X,Y)=Y^4-X^9+X^5Y^2.
\]

One can easily verify that both branches determine the same
Hilbert function.

The branches $(f)$ and $(g)$ have, respectively, the following
Puiseux expansions:
\[
(t^4,t^9-\frac{1}{4}t^{10}-\frac{1}{32}t^{11}+\frac{7}{2048}t^{13}+\frac{1}{512}t^{14}+\frac{39}{65536}t^{15}+\cdots),
\]
and
\[ (t^4,t^9-\frac{1}{4}t^{11}+\frac{1}{32}t^{13}+\frac{1}{128}t^{15}+\cdots),
\]
whose normal forms, which may be obtained through \cite{[Her]},
are, respectively,
\[
(t^4, t^9+t^{10}-\frac{1}{2}t^{11}), \ \ \ \ (t^4, t^9+t^{11}).
\]

According to Table (1), these branches have distinct sets of values of
K\"ahler differentials, which are, respectively,
\[ \{ 4, 8, 9 ,12,13,14,16,\ldots \} \ \ \ {\rm and} \ \ \
\{ 4, 8, 9 ,12,13,15,\ldots \}.
\]

Finally, remark that we have the moduli problem for
branches with multiplicity at most 4 is solved as follows.

For multiplicity 1, 2 and 3, or when $\Lambda\setminus
\Gamma=\emptyset$, this is trivial and already contained in
\cite{[Z2]}.

Suppose now that an equisingularity class is given by a semigroup
$\Gamma$ of multiplicity 4. Then the normal forms in Table (1)
determine a finite family of disjoint constructible sets in affine
spaces, each one corresponding to a set $\Lambda\setminus \Gamma$,
modulo a weighted action of the finite group $G$ of the complex
$(\lambda -v_1)$-th roots of unity, where $\lambda=\min(
\Lambda\setminus \Gamma) -4$.

\quad For example, let $\Gamma=\langle 4, v_1 \rangle$, and
$\Lambda \setminus \Gamma=\{ 2v_1 -4s, \ 3v_1 -4s; \ \ 1\leq s
\leq j-1\}$, for some $j=2,\ldots, \left[ \frac{v_1}{4} \right]$.

So, $\lambda =2v_1-4j$. In this case, the corresponding component
of the moduli is $\mathbb C^{ \ j-1}/G$, where the action of $G$
on $\mathbb C^{ \ j-1}$ is as follows:
$$
(a_{[\frac{v_1}{4}]-j+2},a_{[\frac{v_1}{4}]-j+3}, \ldots,
a_{[\frac{v_1}{4}]}) \sim
(a'_{[\frac{v_1}{4}]-j+2},a'_{[\frac{v_1}{4}]-j+3}, \ldots,
a'_{[\frac{v_1}{4}]}) \ \Leftrightarrow $$
$$
 a'_i=r^{i-v_1}a_i, \ \ \left[\frac{v_1}{4}\right]-j+2\leq  i \leq \left[\frac{v_1}{4}\right], \ r\in G.
$$

Now, if $\Gamma=\langle 4, v_1,v_2\rangle$, then the moduli
consists of one component, which is a point if $v_1=6$, and
$\mathbb C^{ \  [\frac{v_1}{4}]-1}/G$, where $G=\{r\in \mathbb C;
\ \ r^{v_2-2v_1}=1\}$, if $v_1>6$. The action of $G$ is as
follows:
$$
(a_{1},a_{2}, \ldots, a_{[\frac{v_1}{4}]-1}) \sim  (a'_{1},a'_{2},
\ldots, a'_{[\frac{v_1}{4}]-1}) \ \Leftrightarrow
 a'_i=r^{i-v_1}a_i, \ \ 1 \leq  i \leq \left[\frac{v_1}{4}\right]-1, \ r\in G.
$$
\newpage
\begin{center}

{\small \begin{tabular}{|c|l|c|c|} \hline {\bf $\Gamma$} & {\bf
Normal Form} & {\bf $\Lambda\setminus\Gamma$} & $\tau$
\\
\hline $\langle 1 \rangle$ & $\hspace*{2mm} x=t, \ \ \ \ y=t$ &
$\emptyset$  &  \  $0$ \\
\hline $\langle 2,v_1\rangle$ & $\hspace*{2mm} x=t^2, \ \ \ \
y=t^{v_1}$ & $\emptyset$ &  \  $v_1-1$ \\
\hline   &  $\hspace{2mm} x=t^3, \ \ \ \ y=t^{v_1}$ & $\emptyset$
& $2(v_1-1)$
\\ \cline{2-4}
 $\langle 3,v_1\rangle$ &
$\begin{array}{l} x=t^3, \ \ \ \ y=t^{v_1}+t^{2v_1-3j}\\
2\leq j\leq \left [ \frac{v_1}{3}\right ]\end{array}$ & $\begin{array}{l} 2v_1-3s;\\
1\leq s\leq j-1\end{array}$ & $2v_1-j-1$
\\ 
\hline
 &
$ \hspace{2mm} x=t^4, \ \ \ \ y=t^{v_1}$ & $\emptyset$ &
$3(v_1-1)$
\\ \cline{2-4}
 &
$\begin{array}{l} x=t^4,\\
y=t^{v_1}+t^{3v_1-4j}+ \\ \hspace{1cm} \sum_{i=1}^{j-\left [
\frac{v_1}{4}\right ]-2}
{a_it^{2v_1-4\left ( j-\left [ \frac{v_1}{4}\right ] -i\right )}}  \\
2\leq j \leq \left [ \frac{v_1}{2}\right ]
\end{array}$ &
$\begin{array}{c}
3v_1-4s;\\
1\leq s\leq j-1\end{array}$ & $3v_1-j-2$
\\ \cline{2-4}
$\langle 4,v_1\rangle$ & $\begin{array}{l}
x=t^4,\\
y=t^{v_1}+t^{2v_1-4j}+\\ \hspace{1cm}
\sum_{i=k}^{j-2+k}{a_it^{3v_1-4\left ( \left [ \frac{v_1}{4}\right
]+j+1-i\right )}}
\\
2\leq j \leq \left [ \frac{v_1}{4}\right ],\ \ a_k\neq 0, \\
1\leq k\leq \left [ \frac{v_1}{4}\right ]-j
\end{array}$ &
$\begin{array}{l}  2v_1-4s;\\
1\leq s\leq j-1 \\
3v_1-4s;\\
1\leq s\leq\left [ \frac{v_1}{4}\right ] +1-k\end{array}$ &
$3(v_1-1)+k-\left[\frac{v_1}{4}\right]$
\\ \cline{2-4}
 &
$
\begin{array}{l}
x=t^4,\\
y=t^{v_1}+t^{2v_1-4j}+\\ \hspace{1cm} \sum_{i=\left [
\frac{v_1}{4}\right ]-j+1}^{\left [ \frac{v_1}{4}\right ]
-1}{a_it^{3v_1-4\left ( \left [ \frac{v_1}{4}\right ]+j+1-i\right
)}}
\\
2\leq j \leq \left [ \frac{v_1}{4}\right ], \ \ a_{\left [
\frac{v_1}{4}\right ]-j+1}\not =\frac{3v_1-4j}{2v_1}
\end{array}$ &
$\begin{array}{l}
2v_1-4s;\\
1\leq s\leq j-1 \\
3v_1-4s;\\
1\leq s\leq j\end{array}$ & $3v_1-2(j+1)$
\\ \cline{2-4}
 &
$
\begin{array}{l}
x=t^4,\\
y=t^{v_1}+t^{2v_1-4j}+ \frac{3v_1-4j}{2v_1}t^{3v_1-8j}+\\
\hspace{1cm} \sum_{i=\left [ \frac{v_1}{4}\right ]-j+2}^{\left [
\frac{v_1}{4}\right ]}{a_it^{3v_1-4\left ( \left [
\frac{v_1}{4}\right ]+j+1-i\right )}}
\\
2\leq j \leq \left [ \frac{v_1}{4}\right ]
\end{array}$ &
$\begin{array}{c}
2v_1-4s, \ 3v_1-4s;\\
1\leq s\leq j-1
\end{array}$ & $3(v_1-1)-2(j-1)$
\\
\hline $\langle 4,v_1,v_2\rangle$ &
$\begin{array}{l} x=t^4\\
y=t^{v_1}+t^{v_2-v_1}+ \\
\hspace{1cm} \sum_{i=1}^{\left [ \frac{v_1}{4}\right
]-1}{a_it^{v_2-4\left ( \left [ \frac{v_1}{4}\right ] +1-i\right
)}}\end{array}$ & $\begin{array}{c}
v_2+v_1-4s;\\
1\leq s\leq \frac{v_1}{2}-1\end{array}$  & $v_2+\frac{v_1}{2}-2$
\\
\hline
\end{tabular}}
\vspace{3mm}

\noindent {\bf Table (1):} Normal forms for plane branches with
multiplicity less or equal than 4.

\end{center}

\vspace{3mm} The last column of the above table contains the
result of \cite{[NW]}, and is obtained by using the equality
$\tau=c-\sharp \Lambda \setminus \Gamma$, where $c$ is the
conductor of $\Gamma$, which was pointed out in Remark 4.8 of
\cite{[HH1]}.

\end{document}